\documentclass[12pt]{article}
\usepackage[latin9]{inputenc}
\usepackage{amsmath}
\usepackage{amssymb}
\usepackage{esint}

\usepackage{amsthm}\usepackage{amsfonts}\usepackage{latexsym}

\begin{document}

\newtheorem{theo}{Theorem}[section]
\newtheorem{prop}[theo]{Proposition}
\newtheorem{cor}[theo]{Corollary} 
\newtheorem{lem}[theo]{Lemma}
\newtheorem{rem}[theo]{Remark}
\newtheorem{con}[theo]{Conjecture}
\newtheorem{as}[theo]{Assumption}
\newtheorem{de}[theo]{Definition}

\setcounter{page}{1}
\renewcommand{\theequation}{\thesection.\arabic{equation}}

\begin{center}
{\Large \bf 
Phase Transitions \\

in the One-dimensional Coulomb Gas Ensembles\\
}
\end{center}

\begin{center}    
TATYANA S. TUROVA\footnote{
Mathematical Center, University of
Lund, Box 118, Lund S-221 00, 
Sweden.} 
\end{center}

\begin{abstract}

We consider the system of particles on a finite interval with pair-wise nearest neighbours interaction and external force. 
This model was introduced by Malyshev to study the flow of charged particles 
on a rigorous mathematical level. It is a simplified version of a 3-dimensional classical Coulomb gas model.
We study Gibbs distribution at finite positive temperature extending recent  results 
on the  zero temperature case (ground states) with external force. We 
derive the asymptotic for the mean and for the variances
 of the  distances between the neighbouring charges.
We prove that depending on the strength of the external force there are several  phase 
transitions in the local structure of the configuration of the particles in the limit when 
the number of particles goes to infinity.

\end{abstract}

\setcounter{equation}{0}
\renewcommand{\theequation}{\thesection.\arabic{equation}}

\section{Introduction}
Coulomb gas ensembles appears in a variety of  mathematical models. 
Here we focus on the model which describes the 
charges with nearest neighbour Coulomb interaction on an  interval in a presence of external force. This model was introduced and 
studied recently by Malyshev \cite{M}, and then by 
Malyshev and Zamyatin \cite{MZ}. Following the previous authors we call this model 1-dimensional. However, one can view it as 
a simplified version of a 3-dimensional classical 
Coulomb gas model (see, e.g., recent book by Serfaty \cite{Se}), where the particles are assumed to be hardly alined on an interval, imposing only nearest-neighbours interaction.

In another class of related models, which describes, e.g., 
 a one-component plasma (2-dimensional Coulomb gas), 
the charges are related to the
eigenvalues of random normal matrices. Hence, the questions of 
existence and universality of scaling limits for the eigenvalues of a random
normal matrix are in focus of these studies (see, e.g., \cite{A1}, \cite{A2} and
\cite{B}, and the reference therein). Note, that these planar  models (including the 
one-dimensional case) treat only
logarithmic interactions.

The model of \cite{M} and \cite{MZ} is, perhaps, simpler, 
but it is rich enough to exhibit phase transitions which might explain
 certain electric phenomena. We refer to the papers  \cite{M} and \cite{MZ} 
on the relation of this model to physics. 

A local structure of Gibbs configurations 
without external force was analyzed in \cite{MZ}. Note that already in  \cite{M}
it was proved that at the  zero temperature case (ground states) 
there are
 phase transitions in the structure of the configuration of charges under different strength of external force. 
Also for the $d$-dimensional models of Coulomb gas the large deviations principle at speed $n^2$ is proved in \cite{Se} for a rather general situation.

Here we derive the asymptotic mean and the variance of 
the  distances between neighbouring charges. We prove that depending 
on the external force the location of particles undergo phase transition: 
for the weak force (which still may increase with the number of particles) 
the charges remain to be equally spaced over the interval, 
at critical value of the external force they occupy only a finite part of the interval, 
and when the force is above the  critical value all the charges collapse 
in one end of the interval. 

The methods we use here develop the probabilistic  approach of 
\cite{MZ} (and \cite{DF})  but now in an inhomogeneous setting. 

\setcounter{equation}{0}
\renewcommand{\theequation}{\thesection.\arabic{equation}}
\section{Model}

Consider a system of $N+1$ identical particles on the interval $[0,L]$, whose locations 
are represented by a random vector $\bar{Y}=(Y_0,  \ldots , Y_{N})$ where the 
components are ordered: 
\[ 0=Y_{0}<\ldots<Y_{N}=L. \]
Notice here that the lengths of the interval plays no role, so we fix it from now one to be one: $L=1$.
Define a function of potential on the space of vectors 
$\bar{y}=(y_{0}, \ldots , y_{N})$ with $ 0=y_{0}<\ldots<y_{N}=1$:
\begin{equation}\label{U}
U(\bar{y})=  \beta \sum_{k=1}^N
V(y_{k}-y_{k-1})+   \sum_{k=1}^{N}\int_{0}^{y_k}F_sds,
\end{equation}
where positive function $V$ represents a pair-wise  interaction between the particles, 
$\beta>0$ is a parameter, and 
 function $F_s$ represents an external force at point $s\in [0,1]$.

The corresponding 
 Gibbs  distribution of the locations $\bar{Y}$ of particles on the interval is defined by its density
\begin{equation}\label{dG}
f_{\bar{Y}}(\bar{y}):= \frac{1}{Z_{\beta, F}(N)}  e^{-U(\bar{y})} , 
\end{equation}
where  the normalizing factor is 
\begin{equation}\label{Z*}
Z_{\beta, F}(N)=\int\ldots\int_{0
  <
  y_{1}<...<y_{N-1}<
  1}e^{-U(\bar{y})} 
dy_{1}\ldots y_{N-1}.
\end{equation}

\begin{as}\label{A1}
Here we  consider a pair-wise Coulomb repulsive interaction in the form
\begin{equation}
V(x)=\frac{1}{x}, \ \ x>0.  \ \ \ \  \label{V}
\end{equation}

Assume also that the external force $F_y=F$ does not depend on the location $y$,
but it is a function of the number of particles: $F=F(N)$.
 
\end{as}

We study here the asymptotic distribution of the inter-spaces 
between the particles, which are the random variables
$Y_{k}-Y_{k-1}$, $1\leq k\leq N$, when $N \rightarrow \infty$. 

Using Assumptions \ref{A1} let us rewrite the potential function as follows:
\begin{equation}\label{U1}
U(\bar{y})= \beta  \sum_{k=1}^{N} V(y_k-y_{k-1}) + F\sum_{k=1}^{N} (N-k+1)(y_k-y_{k-1})
\end{equation}
\[= \beta  \sum_{k=1}^{N} V(y_{N-k}-y_{N-(k-1)}) + F\sum_{k=1}^{N} k(y_{N-k}-y_{N-(k-1)}). \]
Denote $x_k=y_{N-k}-y_{N-(k-1)}$, $k=1, \ldots, N$. Then 
\[U(\bar{y})= \sum_{k=1}^{N}\left( \beta V(x_k)+ F k x_k \right). \]
This form suggests 
 the following useful representation. 

Let us introduce independent random variables
$X_1, \ldots, X_k$ with density functions 
\begin{equation}\label{fX}
f_{X_k}(x)=\frac{1}{c_k} e^{-\beta V(x)- Fk x} = 
\frac{1}{c_k} e^{-\frac{\beta }{x} - F k x} , \ \ \ x \in [0,1],
\end{equation}
where $c_k$ is the normalizing constant 
\begin{equation}\label{ck}
c_k= \int_0^1  e^{-\frac{\beta }{x} - F k x} dx.
\end{equation}
Then (\ref{U1}) yields
\begin{equation}\label{cond}
X_k \mid _{\sum_{i=1}^{N} X_i =1 } \  \  \stackrel{d}{=}  \  \  Y_{N-k}-Y_{N-(k-1)} .
\end{equation}

\setcounter{equation}{0}
\renewcommand{\theequation}{\thesection.\arabic{equation}}

\section{Results}
Our goal  here is to find the mean and the variance of the conditional distribution (\ref{cond}). The results are listed  in the next theorem in order of increase external force. 
Recall that the model when $F=0$ was treated in \cite{MZ} (it is a particular case of part $(a)$ of the following theorem). 

We shall  use notation 
$g(N)=\Theta(h(N)) $ if  for some constants $0<c<C$
\[ ch(N)<g(N)<Ch(N). \]

\begin{theo}\label{T1}
Consider conditional distribution (\ref{cond}) under Assumptions \ref{A1}.
Define 
\[F_{cr}(N)= 4\beta N.\]
\noindent
(a) If $F(N)=o(N)$ then for all $1\leq k \leq N$
\begin{equation}\label{TE}
\mathbb{E} \left\{ X_k \mid \sum_{i=1}^{N} X_i =1  \right\} = \frac{1}{N} \left(
1-\frac{F}{2\beta N} 
\left(\frac{k}{N} -\frac{1}{2} \right)
 \right)\left(1+ O \left( \frac{F^2}{N^2} \right)   + O \left( \frac{\log N}{\sqrt{N}}\right) \right),
\end{equation}
and
\begin{equation}\label{TVa}
Var\left\{ X_k \mid \sum_{i=1}^{N} X_i =1  \right\} =\frac{1}{2\beta N^3} (1+o(1)).
\end{equation}
\bigskip

\noindent
(b) If $F(N)= F_0 N< F_{cr}(N)$, i.e., 
\begin{equation}\label{cr}
F_0<4 \beta,
\end{equation}
then 
\begin{equation}\label{TE1}
\mathbb{E} \left\{ X_k \mid \sum_{i=1}^{N} X_i =1  \right\} = \frac{1}{a_k N} 
\left(1+ O \left( \frac{\log N}{\sqrt{N}}\right) \right)
\end{equation}
and
\begin{equation}\label{TVb}
Var\left\{ X_k \mid \sum_{i=1}^{N} X_i =1  \right\} =\frac{1}{2\beta a_k^3 N^3} 
\left(1+ O \left( \frac{\log N}{\sqrt{N}}\right) \right),
\end{equation}
where
\begin{equation}\label{ak}
 a_k = \frac{1}{\sqrt{
 1+  \left( \frac{k}{N}-\frac{1}{2}
\right)\frac{F_0}{\beta} + \frac{F_0^2}{16\beta^2 }
} }.
\end{equation}

\bigskip

\noindent
(c) If $F(N)=  F_{cr}(N) =4 \beta N$ 
then 
\begin{equation}\label{crEk}
\mathbb{E} \left\{ X_k \mid \sum_{i=1}^{N} X_i =1  \right\} = 
\sqrt{\frac{1}{4kN+\Theta(N)}}+O(N^{-2/3}),
\end{equation}
and 
\begin{equation}\label{TVc}
Var\left\{ X_k \mid \sum_{i=1}^{N} X_i =1  \right\} = O(N^{-3/2}).
\end{equation}
\medskip

\noindent
(d) If $F(N)= F_0 N> F_{cr}(N)$, i.e., 
\begin{equation}\label{crsup}
F_0> 4 \beta,
\end{equation}
then 
\begin{equation}\label{cEXsup}
\mathbb{E} \{ X_1 \mid \sum_{i=1}^{N} X_i =1 \} 
= 1-\sqrt{\frac{4 \beta}{F_0} }+O\left(\frac{1}{\sqrt{N}}\right),
\end{equation}
\begin{equation}\label{supV1}
Var\left\{ X_1 \mid \sum_{i=1}^{N} X_i =1 \right\} 
= O\left(\frac{1}{N^{3/4}}\right),
\end{equation}
and for all $k\geq 2$
\begin{equation}\label{supE}
\mathbb{E} \{ X_k\mid \sum_{i=1}^{N} X_i =1 \} 
= \sqrt{\frac{\beta}{(k-1)F_0N +\lambda_0}}
+O
\left( \frac{\log N}{ (kN)^{3/4} } \right), 
\end{equation}
and 
\begin{equation}\label{supV}
Var\left\{ X_k \mid \sum_{i=1}^{N} X_i =1 \right\}  =O\left( \frac{1}{(kN)^{3/2}} \right),
\end{equation}
where $\lambda_0= \lambda_0(\beta, F_0)$ is the unique solution to 
\begin{equation}\label{M24}
\frac{\int_0^1 xe^{\lambda_0 x -\frac{\beta}{x}}dx}{\int_0^1 e^{\lambda_0 x -\frac{\beta}{x}}dx}=1-\sqrt{\frac{4 \beta}{F_0} }.
\end{equation}
\medskip

\noindent
(e) If $F(N)\gg  N$
then 
\begin{equation}\label{supp1}
\mathbb{E} \left\{ X_1 \mid \sum_{i=1}^{N} X_i =1  \right\} = 
 1-\sqrt{\frac{4 \beta N}{F} }+O(F^{-2/3}),
\end{equation}
and 
\begin{equation}\label{supp2}
Var\left\{ X_1 \mid \sum_{i=1}^{N} X_i =1  \right\} = O(F^{-3/2}).
\end{equation}

\end{theo}

This theorem confirms that the phase transitions in the ground states discovered in \cite{M} take place for the Gibbs measure at any positive finite temperature as well. Furthermore, here we find also that even the fluctuations (which are not present of course at ground states) around the mean, i.e., the scalings of the variances 
also undergo phase transition. 

Theorem \ref{T1} describes 5 qualitatively different states for the locations of particlaes in the limit. These are the following. 

1. In case $(a)$ (weak subcritical force) the particles remain to be equally spaced on the average 
at intervals of equal length $N^{-1}$, just as if $F=0$, 
and the variances of the spacings are minimal, they are of order $N^{-3}$.

2. In case $(b)$ (medium subcritical force) the same orders  $N^{-1}$ and 
 $N^{-3}$ for the mean and the variance of the spacings are preserved, however, not homogeneous any longer. 
The constants $a_k$  take different values between 
$a_{1}= \frac{1}{1-\frac{F_0}{4\beta}}(1+o(1/N))>1$ 
and $\frac{1}{2}<a_N= \frac{1}{1+\frac{F_0}{4\beta}}<1$.
In particular, when $F_0 \uparrow 4\beta$ one has $a_N \rightarrow 1/2$, while 
 $a_1 \rightarrow \infty$.

3. In case $(c)$ (critical force) all 
the spacings converge to zero (in $L^2$ at least), hence, the particles still densly cover the entire interval, 
however the order of the mean varies monotonically 
between $N^{-1/2}$ to $N^{-1}$.

4. In case $(d)$ (weak supercritical force) the particles occupy densely only a positive fraction of the interval, while on the remaining fraction there are no particles (except the ones with  fixed positions at the ends). 

5. In the case $(e)$ (strong supercritical force) all the particles 
(except the ones with  fixed positions at the ends) are jammed towards one end, the phenomena reminding a condensation.

\setcounter{equation}{0}
\renewcommand{\theequation}{\thesection.\arabic{equation}}
\section{Proof of Theorem \ref{T1}.} 

Consider the conditional density of $X_k$. Following ideas of \cite{DF} we shall introduce another parameter into the distribution of $X_k$. Namely, for any 
$\lambda \in \mathbb{R}$ define a density 

\begin{equation}\label{fXl}
f_{k,\lambda}(u)=
\frac{1}{c_k(\lambda)} e^{-\frac{\beta }{u} - (\lambda +F k) u } , \ \ \ u  \in [0,1] ,
\end{equation}
\[c_k(\lambda)= \int_0^1 e^{-\frac{\beta }{u} - (\lambda +F k) u } du,\]
and denote the corresponding random variable $X_{k,\lambda}$. In these notations 
$X_k\stackrel{d}{=}  
X_{k,0}$.
We assume, that for each fixed 
$\lambda$ the random variables 
$X_{k,\lambda}, k=1, \ldots, N,$ are independent. 
The remarkable property of these random variables is that 
for any $\lambda \in \mathbb{R}$
the following equality in distribution holds:
\begin{equation}\label{condl}
X_k \mid _{\sum_{i=1}^{N} X_i =1 } \  \  \stackrel{d}{=}  \  \  X_{k,\lambda} \mid _{\sum_{i=1}^{N} X_{i,\lambda} =1 }.
\end{equation}
Indeed, 
denoting
\[S_N=\sum_{i=1}^NX_i , \ \ \ S_{N, {\hat{k}}}=\sum_{i\neq k}X_i\]
and, correspondingly,
\[S_{N, \lambda} =\sum_{i=1}^NX_{i,\lambda}, \ \ \ S_{N, {\hat{k}}, \lambda}=\sum_{i\neq k}X_{i,\lambda},\]
it is straightforward to check that for the conditional densities one has
\begin{equation}\label{cfX}
f_{X_k \mid {\sum_{i=1}^{N} X_i =1 }} (x):= 
\frac{f_{X_k}(x) f_{S_{N, {\hat{k}}}} (1-x)}{f_{S_{N}} (1)}
 \  \  =  \  \  \frac{f_{k,\lambda}(x) f_{S_{N, {\hat{k}}, \lambda}} (1-x)}{f_{S_{N,\lambda}} (1)},
\end{equation}
i.e., the right-hand side does not depend on $\lambda$.

Notice, that Diaconis and Freedman \cite{DF} used this property for the identically distributed random variables, their result was then used in \cite{MZ} to treat the case without external force. Here we show that a similar argument works even without condition on the identity of the distributions. The main idea is to tune the free parameter $\lambda$ so that the condition $
\sum_{k=1}^N X_{k,\lambda}
=1$ will be natural in the following sense. Since we have a sum of independent random variables,
we may expect that 
due to the Central Limit Theorem  the density of the 
normalized sum  $S_{N,\lambda}$ converges to the  density of a normal distribution with the maximum at the point of its expected value. Therefore we shall  choose $\lambda=\lambda(N)$ so that 
\begin{equation}\label{ES}
\mathbb{E} S_{N, \lambda} =\sum_{k=1}^N\mathbb{E} X_{k,\lambda} = 1.
\end{equation}
Notice that  the definition (\ref{fXl}) yields (see also \cite{MZ} and \cite{DF}) that 
 $\mathbb{E} X_{k,\lambda}$ is a strictly decreasing function of $\lambda$, such that 
$\mathbb{E} X_{k,\lambda} \rightarrow 0$ when $\lambda \rightarrow \infty$ while $\mathbb{E} X_{k,\lambda} \rightarrow 1$ when $\lambda \rightarrow - \infty$. Hence, equation (\ref{ES}) defines  uniquely  $\lambda=\lambda(N)$ which satisfies this condition. To solve equation 
(\ref{ES}) first we consider $\mathbb{E} X_{k,\lambda}$. Notice, that in \cite{MZ} one can find  the principal term of the asymptotic of this value. Here using the  arguments of \cite{MZ} 
 we get more details.

\begin{prop}\label{P1} For any  $1\leq k \leq N$ and $F\geq 0$
\begin{equation}\label{mk}
m_{k,\lambda} := \mathbb{E} X_{k,\lambda} = 
\sqrt{\frac{\beta}{kF+\lambda} }+ \frac{3}{4(kF+\lambda)} (1+o(1)) ,
\end{equation}
and 
\begin{equation}\label{M16}
Var (X_{k,\lambda})=\frac{\sqrt{\beta}}{2(kF+\lambda)^{3/2}}(1+o(1)).
\end{equation}
as $kF+\lambda \rightarrow \infty$.

\end{prop}

\noindent
{\bf Proof.}
Let us write here
\begin{equation}\label{M2}
 \lambda_k = \lambda +kF,
\end{equation}
and, correspondingly, $X_{k,\lambda}=X_{\lambda_k}$ (this notations are  consistent with the ones in \cite{MZ}).

For any  $\lambda >0 $ and $\alpha\in \{1,2,3\}$ define 
\begin{equation}\label{Ia}
I_{\alpha}(\lambda, \beta)= \int_0^{\infty}x^{\alpha-1}e^{-\lambda x -\frac{\beta }{x}}dx.
\end{equation}
With this notation we have
\begin{equation}\label{M1}
\mathbb{E} X_{\lambda}=\frac{\int_0^1 x e^{-\lambda x -
\frac{\beta}{x}} dx}{\int_0^1  e^{-\lambda x -\frac{\beta}{x}} dx}=
\frac{
I_{2}(\lambda, \beta) +O(e^{-\frac{1}{2}\lambda})
}{I_{1}(\lambda, \beta) +O(e^{-\frac{1}{2}\lambda}) } .
\end{equation}
It is shown in \cite{MZ} that 
\begin{equation}\label{M13}
I_{\alpha}(\lambda, \beta)= \frac{2\beta^{\alpha/2} K_{\alpha}(2\sqrt{\lambda\beta})}{\lambda^{\alpha/2}},
\end{equation}
where $K_{\alpha}(z)$ is a Bessel function for which the asymptotic expansion when $z\rightarrow \infty$ is known to be
\begin{equation}\label{M12}
K_{\alpha}(z) = \sqrt{\frac{\pi}{2}}\frac{e^{-z}}{\sqrt{z}}\left( 1+ \frac{4\alpha^2 -1}{8z} +o(z^{-1})\right).
\end{equation}
Hence, using (\ref{M13}) and then (\ref{M12}), we derive from  (\ref{M1}) 
\begin{equation}\label{M14}
\mathbb{E} X_{\lambda}=
\frac{\frac{\beta}{\lambda}
K_{2}(2\sqrt{\lambda\beta}) +O(e^{-\frac{1}{2}\lambda})
}{ \sqrt{\frac{\beta}{\lambda}}
K_{1}(2\sqrt{\lambda\beta}) +O(e^{-\frac{1}{2}\lambda}) } 
\end{equation}
\[
=\frac{ \frac{\beta}{\lambda}
\left( 1+ \frac{15}{16\sqrt{\lambda\beta}} +o(\lambda^{-1/2}) \right)
 +O(\lambda^{1/4}e^{-\frac{1}{2}\lambda+2\sqrt{\lambda\beta}})
}{ \sqrt{\frac{\beta}{\lambda}}
 \left( 1+ \frac{3}{16\sqrt{\lambda\beta}} +o(\lambda^{-1/2}) \right)
 +O(\lambda^{1/4}e^{-\frac{1}{2}\lambda+2\sqrt{\lambda\beta}})
 } =
\sqrt{\frac{\beta}{\lambda}}\ \
\frac{ 
 1+ \frac{15}{16\sqrt{\lambda\beta}} +o(\lambda^{-1/2}) 
}{ 
  1+ \frac{3}{16\sqrt{\lambda\beta}} +o(\lambda^{-1/2}) 
 } 
\]
\[
=\sqrt{\frac{\beta}{\lambda} }\left( 1+ \frac{3}{4\sqrt{\lambda\beta}} +o(\lambda^{-1/2}) \right)=\sqrt{\frac{\beta}{\lambda} }+\frac{3}{4\lambda} +o(\lambda^{-1}),
\]
which yields (\ref{mk}).

In a similar manner we derive as well
\begin{equation}\label{M15}
\mathbb{E} (X_{\lambda})^2=\frac{\int_0^1 x^2 e^{-\lambda x -
\frac{\beta}{x}} dx}{\int_0^1  e^{-\lambda x -\frac{\beta}{x}} dx}=
\frac{
I_{3}(\lambda) +O(e^{-\frac{1}{2}\lambda})
}{I_{1}(\lambda) +O(e^{-\frac{1}{2}\lambda}) } 
\end{equation}
\[=
\frac{\beta}{\lambda}\ \
\frac{ 
 1+ \frac{35}{16\sqrt{\lambda\beta}} +o(\lambda^{-1/2}) 
}{ 
  1+ \frac{3}{16\sqrt{\lambda\beta}} +o(\lambda^{-1/2}) 
 } = \frac{\beta}{\lambda}\left( 1+ \frac{2}{\sqrt{\lambda\beta}} +o(\lambda^{-1/2}) \right),
\]
which together with (\ref{M14}) yields as well
\[
Var (X_{\lambda})=\frac{\sqrt{\beta}}{2\lambda^{3/2}}(1+o(1)),
\]
and (\ref{M16}) follows. Proposition is proved. \hfill$\Box$

Now we can choose $\lambda$ so that (\ref{ES}) holds. We shall describe the asymptotic of this value. 

\begin{lem}\label{L1} 
Assume that   $\lambda=\lambda(N,F)$ is chosen so that  (\ref{ES}) holds, i.e., 
\[\sum_{k=1}^N \mathbb{E} X_{k,\lambda} =1.\]

\noindent
(a) If  $F=o(N)$ then
\begin{equation}\label{la1}
\lambda(N,F) =\beta \left( 1-\frac{F}{4\beta N}
\right)^2 N^2  \left( 1+\frac{1}{ N}
\right).
\end{equation}

\noindent
(b) If  $F=F_0 N$ and $F_0<4\beta$ then 
\begin{equation}\label{la2}
\lambda(N,F) =\beta \left( 1-\frac{F}{4\beta N}
\right)^2 N^2+ \Theta \left(\frac{N}{\log N}
\right).
\end{equation}

\noindent
(c)  If  $F= 4\beta N$ then 
\begin{equation}\label{la2c}
\lambda(N,F) = \Theta \left( 
N
\right)>0.
\end{equation}

\noindent
(d) If  $F=F_0 N$ where 
$F_0>4\beta$, including  $F_0=F_o(N)\gg 1$, then 
\begin{equation}\label{la3}
\lambda(N,F) = - F_0 N + \lambda_0+o(1),
\end{equation}
where $\lambda_0= \lambda_0(\beta, F_0)$ is the unique solution to (\ref{M24}).
\end{lem}
 
\begin{rem}\label{R1}
In the case $(d)$  function $\lambda_0= \lambda_0(\beta, F_0)$ is increasing in 
$F_0$; in particular, $\lambda_0 (F_0) = \sqrt{\frac{F_0}{4\beta}}(1+o(1)) \rightarrow \infty$ if $F_0 \rightarrow \infty$.
\end{rem}

\noindent
{\bf Proof.}
Let us solve (\ref{ES}) when $F> 0$. First 
using formula (\ref{mk})  we derive
\begin{equation}\label{mkS}
\sum_{k=1}^N \mathbb{E} X_{k,\lambda} = 
\sum_{k=1}^N \left(\sqrt{\frac{\beta}{kF+\lambda} }+ \frac{3}{4(kF+\lambda)} (1+o(1)) \right).
\end{equation}
Consider now separately different cases.

If $F=o(N)$ we derive from equations (\ref{mkS}) and  (\ref{ES})
 \begin{equation}\label{mk1*}
1=\sqrt{\frac{\beta}{F} }\left( \int_1^N \frac{1}{\sqrt{x+\frac{\lambda}{F}}}dx + O\left(
\frac{1}{\sqrt{1+\frac{\lambda}{F}}}
\right)\right)
+ O\left(
\frac{N}{\lambda}
\right)
\end{equation}
\[ = 2 \sqrt{\beta}\ \frac{N-1}{\sqrt{FN+\lambda}
+\sqrt{F+\lambda}}  
 + O\left(
\frac{1}{\sqrt{F+{\lambda}}}
\right) + O\left(
\frac{N}{\lambda}
\right)\]
\[ = 2 \sqrt{\beta} \ \frac{N}{\sqrt{FN+\lambda}
+\sqrt{F+\lambda}} \left(1
 + O\left(
\frac{1}{N}
\right)\right).\]
This gives us
\begin{equation}\label{mk2}
 \lambda(N,F)=\beta N^2 \left(1- \frac{F}{4\beta N}\right)^2 (1+O(1/N)),
\end{equation}
where the last term is uniform in $F$. This confirms statement $(a)$.

Assume now that  $F(N)=F_0N$ where $F_0$  is some positive constant. 
Then we derive from equation (\ref{mkS}) 
 \begin{equation}\label{mk1}
\sum_{k=1}^N \mathbb{E} X_{k,\lambda} =\sqrt{\frac{\beta}{F_0N} }\left( \int_1^N \frac{1}{\sqrt{x+\frac{\lambda}{F_0N}}}dx + \Theta \left(
\frac{1}{\sqrt{1+\frac{\lambda}{F_0N}}}
\right)\right)
+ \Theta \left(
\frac{\log N}{N}
\right).
\end{equation}
In this case equation (\ref{ES}) is equivalent to 
\begin{equation}\label{M20}
\sqrt{\frac{F_0}{4 \beta} }
= 
\sqrt{1+\frac{\lambda}{F_0N^2} }
-\sqrt{ \frac{1}{N} +\frac{\lambda}{F_0N^2} }
 + \Theta \left(
\frac{1}{\sqrt{N+\frac{\lambda}{F_0}}}
\right)
+ \Theta \left(
\frac{\log N}{{N}}
\right) .
\end{equation}

Notice that the function on the right in  (\ref{M20}) is decreasing in $\lambda$, and therefore the principal term on the left is at most 1. 
(In  particular, this yields that 
 if $F_0>4\beta$ there is no positive solution $\lambda= \Theta \left(N^2\right)$.)

Set 
\begin{equation}\label{J0}
\lambda=xF_0N^2
\end{equation}
 in  (\ref{M20}). It is straightforward to  derive
\begin{equation}\label{J1}
x=\frac{B^2-\frac{4}{N}}{4(1+\frac{1}{N}-B)},
\end{equation}
where
\[B= -\left(\sqrt{\frac{F_0}{4\beta} }+\Theta\left(
\frac{1}{\sqrt{N+xN^2}}
\right)
+ \Theta \left(
\frac{\log N}{{N}}
\right) \right)^2+1+\frac{1}{N}\]
\[  = 1-\frac{F_0}{4\beta} +\Theta\left(
\frac{1}{\sqrt{N+xN^2}}
\right)
+ \Theta \left(
\frac{\log N}{{N}}
\right) . \]
Hence, if $F_0<4\beta$
\begin{equation}\label{J2}
x=\frac{\beta}{F_0}\left(1+\Theta\left(
\frac{1}{\sqrt{N+xN^2}}
\right)
+ \Theta \left(
\frac{\log N}{{N}}
\right)  \right)
\left(1-\frac{F_0}{4\beta}
\right)^2, 
\end{equation}
which together with (\ref{J0}) proves statement (b).

If $F_0=4\beta$, then  (\ref{J1}) gives us
\[x=\frac{\beta}{F_0}\left(1+\Theta\left(
\frac{1}{\sqrt{N+xN^2}}
\right)
+ \Theta \left(
\frac{\log N}{{N}}
\right)  \right)
\left(\Theta\left(
\frac{1}{\sqrt{N+xN^2}}
\right)
+ \Theta \left(
\frac{\log N}{{N}}
\right) 
\right)^2\]
\[=\Theta\left(
\frac{1}{{N+xN^2}}
\right)
+ \Theta \left(
\frac{\log N}{{N}}
\right)^2 ,\]
which yields
\begin{equation}\label{J3}
x= \Theta \left(
\frac{1}{{N}}
\right) ,
\end{equation}
and notably this is a positive function. Hence, statement (c) of the Lemma \ref{L1}
follows by (\ref{J0}) and (\ref{J3}).

Let  $F_0>4\beta$ (here $F_0$ can be a constant or even an increasing function of $N$).
First we observe that if $F_0>4\beta$ and $\lambda>0$ then by  (\ref{mk1}) 
(with a help of (\ref{M20}))
we have $\sum_{k=1}^N \mathbb{E} X_{k,\lambda} <1.$ 
Therefore the (unique) solution $\lambda$ to the equation
\[\sum_{k=1}^N \mathbb{E} X_{k,\lambda} =1\]
is negative in this case, we shall derive it now. 

Consider again formula (\ref{mkS}) with $F=F_0 N$:
\begin{equation}\label{ES2}
\sum_{k=1}^N \mathbb{E} X_{k,\lambda} = \mathbb{E} X_{1,\lambda} + 
\sum_{k=2}^N \left(\sqrt{\frac{\beta}{kF+\lambda} }+ \frac{3}{4(kF+\lambda)} (1+o(1)) \right),
\end{equation}
where similar to  (\ref{M20}) and (\ref{mk1}) for any $\lambda>-F_0 N$ and such that $|\lambda |\leq 2 F_0 N $
\[
\sum_{k=2}^N \mathbb{E} X_{k,\lambda}
= \sqrt{\frac{4 \beta}{F_0} }
\left( \sqrt{1+\frac{\lambda}{F_0N^2} }
-\sqrt{ \frac{2}{N} +\frac{\lambda}{F_0N^2} }\right)
 + O\left(
\frac{1}{\sqrt{2F_0 N+\lambda}}
\right)
+ O\left(
\frac{\log N}{{N}}
\right)\]
\begin{equation}\label{M21}
\rightarrow \sqrt{\frac{4 \beta}{F_0} }<1 \ \ \mbox{ as } N \rightarrow \infty.
\end{equation}
Let us find a  (negative) $\lambda$ which satisfies 
\begin{equation}\label{M22}
\mathbb{E} X_{1,\lambda} =1-\sqrt{\frac{4 \beta}{F_0} }.
\end{equation}
Set now $\lambda=-F_0 N + \lambda_0$, where 
$\lambda_0$ satisfies (\ref{M24}). Then by the definition (\ref{fXl}) 
\[\mathbb{E} X_{1,\lambda} = \frac{\int_0^1 xe^{-\lambda_0 x -\frac{\beta}{x}}dx}{\int_0^1 e^{-\lambda_0 x -\frac{\beta}{x}}dx} =1-\sqrt{\frac{4 \beta}{F_0} },\]
which yields the  desired property (\ref{M22}). Observe also that
(\ref{M21}) holds for $\lambda=-F_0 N + \lambda_0$. This yields  statement $(d)$
of the lemma, where the correction term $o(1)$ is due to the asymptotically convergence in (\ref{M21}) of the sum $\sum_{k=2}^N \mathbb{E} X_{k,\lambda}.$ Lemma is proved. \hfill$\Box$ 

\begin{rem}\label{R2} 
It follows from (\ref{M22}) that if $F_0=F_0(N)$ is unbounded, i.e., when $F\gg N$,
for the chosen $\lambda$ we have $X_{1,\lambda} \rightarrow 1$, while $\sum_{k=2}^N \mathbb{E} X_{k,\lambda} \rightarrow 0
$ as $N\rightarrow \infty$.
\end{rem}

\begin{cor}\label{C1} 
Assume 
  $\lambda=\lambda(N,F)$ satisfies (\ref{ES}) and is chosen as in Lemma \ref{L1}.  

$(\mbox{I})$  If $\lim _{N\rightarrow \infty} F/N = F_0 <4\beta$ then there is a positive constant $C=C(\beta, F_0)$ such that 
\begin{equation}\label{varS}
\sigma_N^2=\sum_{k=1}^N Var (X_{k,\lambda})=
\frac{C}{  N^2 }(1+o(1)),
\end{equation}
where for any $k$
\begin{equation}\label{vark}
Var (X_{k,\lambda})= \frac{\sqrt{\beta}}{2\left(kF+
\beta\left(1-\frac{F_0}{4\beta} \right)^2N^2
\right)^{3/2}}(1+o(1))=
O(N^{-3}).
\end{equation}

\bigskip
$(\mbox{II})$ 
If $\lim _{N\rightarrow \infty} F/N =4\beta $ then 
\begin{equation}\label{Var1a}
\sum_{k=1}^N Var (X_{k,\lambda})= \Theta \left( N^{-3/2}\right),
\end{equation}
where for any $k$
\begin{equation}\label{vark2}
Var (X_{k,\lambda})= 
\frac{\sqrt{\beta}}{2(kF+\Theta\left(N\right))^{3/2}}(1+o(1)).
\end{equation}

\bigskip
$(\mbox{III})$ 
If $\lim _{N\rightarrow \infty} F/N >  4\beta$ then there is a positive constant $C=C(\beta)$ such that 
\begin{equation}\label{varS1}
\sum_{k=2}^N Var (X_{k,\lambda})=
\frac{C}{  F^{3/2} }(1+o(1))=O(N^{-3/2}), 
\end{equation}
where
  for any $k\geq 2$
\begin{equation}\label{Vark}
Var (X_{k,\lambda})= \frac{\sqrt{\beta}}{ 2 (k-1)^{3/2}F^{3/2}}(1+o(1)),
\end{equation}
while

$(a)$ if $\lim _{N\rightarrow \infty} F/N =F_0$ and $ F_0>4\beta $  is a positive constant, then 
\begin{equation}\label{Var1b}
Var (X_{1,\lambda})=\sigma(F_0) 
\end{equation}
 is also some positive constant; 

$(b)$
if 
$F\gg N$ then 
\begin{equation}\label{Var1c}
Var (X_{1,\lambda})=\frac{\sqrt{2}\beta^{5/4}}{F^{3/4}}(1+o(1)).
\end{equation}
\end{cor}

\noindent
{\bf Proof.}
By  (\ref{M16}) we have for $\lambda=\lambda(N,F)$
\[Var (X_{k,\lambda})=\frac{\sqrt{\beta}}{2(kF+\lambda)^{3/2}}(1+o(1)).\]
In case $(\mbox{I})$ we have $\lambda= \Theta(N^2)$, thus (\ref{Vark}) follows.
Then we derive 
using (\ref{M16}) we derive
\begin{equation}\label{M19}
\sum_{k=1}^N Var (X_{k,\lambda})=\sum_{k=1}^N \frac{\sqrt{\beta}}{2(kF+\lambda)^{3/2}}(1+o(1))= 
\frac{N}{\lambda^{3/2}} \ 
\frac{\sqrt{\beta} }{  \frac{NF}{\lambda} +1 + \sqrt{\frac{NF}{\lambda} +1}
}(1+o(1)) 
\end{equation}
as $\lambda \rightarrow \infty$. 
Then the statement follows by Lemma \ref{L1}. 

Similarly one treats  the remaining cases, taking  into account Lemma \ref{L1}
and also Remark \ref{R1}. \hfill$\Box$

\bigskip

Consider $$S_{N,\lambda} = \sum_{k=1}^N X_{k, \lambda}.$$
From now on we assume that $\lambda=\lambda (N,F)$ is chosen so that $\mathbb{E}S_{N,\lambda} =1 $ (see Lemma \ref{L1}.) We shall also use notations
\[X_{k, \lambda}= X_{\lambda_k},\]
where as above $\lambda_k=kF+\lambda (N,F)$.

Define a random variable
\begin{equation}\label{ZN}
Z_N=\frac{S_{N,\lambda}-1}{\sqrt{Var(S_{N,\lambda})}}=
\frac{1}{{\sigma_N}}\sum_{k=1}^N Y_k,
\end{equation}
where $Y_k = X_{k,\lambda}-m_{k,\lambda}$,
\[m_{k,\lambda}= \mathbb{E}X_{k,\lambda}, \ \ \sigma_N = 
\sqrt{Var(S_{N,\lambda})} = \sqrt{ \sum_{k=1}^N  
Var(X_{k,\lambda})
}.\]
We denote  $\phi_{\cdot}$ and $f_{\cdot}$, correspondingly
the characteristic function and the  density for the random variables. 
Following method  of \cite{MZ} we shall prove that $f_{Z_N}$ converges to the normal density. The only difference is that here we are dealing with non-identically distributed random variables. 

\begin{lem}\label{L2}
Assume that $\lim _{N\rightarrow \infty} F/N = F_0<4\beta $, and let
  $\lambda=\lambda(N,F)$ satisfy (\ref{ES}).  Then 
\begin{equation}\label{M47}
\left|f_{Z_N}(x)-\frac{1}{\sqrt{2\pi}}e^{-x^2/2}\right|\leq O(N^{-1/2}).
\end{equation}
\end{lem}

\noindent
{\bf  Proof.}
Consider the following uniform in $x$ bound (follows by the 
Fourier inverse formula for the densities)
\begin{equation}\label{M23}
|f_{Z_N}(x)-\frac{1}{\sqrt{2\pi}}e^{-x^2/2}|\leq \frac{1}{\sqrt{2\pi} }
\int_{-\infty}^{\infty}|\phi_{Z_N}(t)  -e^{-t^2/2}|dt.
\end{equation}
We shall use the following result of Petrov \cite{P}. 
\medskip

\noindent
{\bf Lemma} (\cite{P}, p.109) 
{\it Let
\[L_N:=\frac{\sum_{k=1}^N \mathbb{E} |Y_k|^3}{{\sigma^3_N}}.\]
Then for all } $|t|\leq \frac{1}{4 L_N}$ 
\begin{equation}\label{Pe}
|\phi_{Z_N}(t)   -e^{-t^2/2}|\leq 16 L_N |t^3|e^{-t^2/3}.
\end{equation}
\medskip

To make use of the bound (\ref{Pe}) we 
have to show that $L_N$ a decreasing towards zero as $N\rightarrow \infty$. 

\begin{prop}\label{P2}
Under assumptions of  Lemma \ref{L2}  one has 
\begin{equation}\label{M32}
L_N=\frac{\sum_{k=1}^N \mathbb{E} |Y_k|^3}{{\sigma^3_N}}=O\left(\frac{1}{\sqrt{N}}
\right).
\end{equation}
\end{prop}

\noindent
{\bf Proof.}
Consider first for any fixed $k$
\begin{equation}\label{M25}
\mathbb{E} |Y_k|^3=\mathbb{E}|X_{k,\lambda}-m_{k,\lambda}|^3=\frac{\int_0^1 |x-m_{k,\lambda}|^3 e^{-\lambda_k x -\frac{\beta}{x}}dx}{\int_0^1 e^{-\lambda_k x -\frac{\beta}{x}}dx},
\end{equation}
where  $\lambda_k=\lambda +kF$, and by Lemma \ref{L1} we have here 
$\lambda_k\rightarrow \infty$. Recall that by Proposition \ref{P1} 
\begin{equation}\label{M27}
m_{k,\lambda}=\sqrt{\frac{\beta}{\lambda_k}}+O\left( \frac{\beta}{\lambda_k}\right)
\end{equation}
when $\lambda_k\rightarrow \infty$. 
Let us write here $m_{k,\lambda}=m(\lambda_k)$, and consider 
\begin{equation}\label{M26}
J(\lambda):
=\frac{\int_0^1 |x-m(\lambda)|^3 e^{-\lambda x -\frac{\beta}{x}}dx}{\int_0^1 e^{-\lambda x -\frac{\beta}{x}}dx}
\end{equation}
for large  $\lambda$. Let us define a function
\begin{equation}\label{S}
s(x)=\lambda x + \frac{\beta}{x}, \ \ x>0.
\end{equation}
Denote $x_{0}$ the argument of the minimal value of $s(x)$ for $x>0$:
\begin{equation}\label{M3}
x_0= \sqrt{\frac{ \beta }{\lambda}}, 
\end{equation}
where 
\begin{equation}\label{M4}
s'(x_{0}) = 0, \ \ \mbox{and } \ \ 
s''(x_{0}) = \frac{ 2\beta }{x_{0}^3} = 2\frac{\lambda^{3/2}}{\sqrt{\beta} } .
\end{equation}
It is straightforward to compute that
 for any $\varepsilon>-\sqrt{\lambda}x_{0} $ and for all large $\lambda$
\begin{equation}\label{M5}
s\left(x_{0}+\frac{\varepsilon}{\sqrt{\lambda}}\right) 
\geq  s(x_{0}) + 
\frac{\sqrt{\lambda}\varepsilon^2 }{2(\sqrt{\beta} +|\varepsilon|)}.
\end{equation}
Since $s'(x)<0$ if $x<x_{0}$ and $s'(x)>0$ if $x>x_{0}$, the bounds (\ref{M5}) 
and  (\ref{M27}) 
imply
for any $0 \leq \alpha \leq 3$
\begin{equation}\label{M6}
\int_0^1 |x-m(\lambda)|^{\alpha} e^{-s(x)}dx =
\int_{
x_{0}- \frac{\varepsilon}{\sqrt{\lambda}} 
}^{
x_{0}+\frac{\varepsilon}{\sqrt{\lambda}}
}  |x-x_0|^{\alpha} e^{-s (x)}dx 
\end{equation}
\[+\alpha \ O\left(\left( \frac{\varepsilon}{\sqrt{\lambda}} \right)^{\alpha} \frac{1}{\lambda}\right)
\int_0^1  e^{-s(x)}dx
+ e^{-s (x_{0} )}O\left(e^{- \frac{\sqrt{\lambda}\varepsilon^2 }{2(\sqrt{\beta} +|\varepsilon|)}} \right).
\]
Next for $\varepsilon=o(1)$ we derive using formulas (\ref{M4}) and 
(\ref{M3}) 
\begin{equation}\label{M7}
\int_{
x_{0}- \frac{\varepsilon}{\sqrt{\lambda}} 
}^{
x_{0}+\frac{\varepsilon}{\sqrt{\lambda}}
} |x-x_0|^{\alpha} e^{-s (x)}dx = e^{-s(x_{0} )} \int_{
x_{0}- \frac{\varepsilon}{\sqrt{\lambda}} 
}^{
x_{0}+\frac{\varepsilon}{\sqrt{\lambda}}
}|x-x_0|^{\alpha} e^{-\frac{1}{2}s'' (x_{0} )(x-x_{0})^2 \left(1+O\left(\frac{x-x_{0}}{x_{0}}\right) \right) }dx
\end{equation}
\[= e^{-s (x_{0} )}(1+O( \varepsilon)) \left(\frac{1}{\sqrt{s'' (x_{0})} } \right)^{\alpha +1}
\int_{
- \frac{\varepsilon}{\sqrt{\lambda}} \sqrt{ s'' (x_{0} )(1+O( \varepsilon)) } 
}^{\frac{\varepsilon}{\sqrt{\lambda}}  \sqrt{ s'' (x_{0} )(1+O( \varepsilon))} 
}\ |x|^{\alpha} e^{-\frac{1}{2}x^2} dx \]

\[=  e^{-s (x_{0} )}  (1+O( \varepsilon))\left(\frac{1}{\sqrt{s'' (x_{0})} } \right)^{\alpha +1}
 \left(1+O\left(
e^{-\frac{\varepsilon^2 s'' (x_{0} )}{4\lambda}} 
\right) \right) \int_{-\infty}^{\infty}|x|^{\alpha} e^{-\frac{1}{2}x^2} dx
.\]
We can choose now $\varepsilon=(\log{\lambda})^2 /\lambda^{1/4}$  so that 
 (\ref{M7}) combined with  (\ref{M6}) gives us 
\begin{equation}\label{M9}
\int_0^1 |x-m(\lambda)|^{\alpha} e^{-s(x)}dx =
e^{-s (x_{0} )}  (1+O( \varepsilon))\left(\frac{1}{\sqrt{s'' (x_{0})} } \right)^{\alpha +1}
\int_{-\infty}^{\infty}|x|^{\alpha} e^{-\frac{1}{2}x^2} dx
\end{equation}
\[+\alpha O\left( 
\frac{   1  }{   \lambda^{1+\frac{\alpha}{2} }   } 
\right)
\int_0^1  e^{-s(x)}dx.\]

Making use of the last formula with $\alpha=3$ and $\alpha=0$  in 
 (\ref{M26}),
and taking into account (\ref{M4})  we derive for all  $\lambda$ 
\begin{equation}\label{M34}
J(\lambda)=c (1+o(1))\left(\frac{1}{\sqrt{s'' (x_{0})} } \right)^{3 }=c (1+o(1))
 \left(\frac{\sqrt{\beta} }{2\lambda^{3/2}} \right)^{3/2 } + O\left( 
\frac{   1  }{   \lambda^{1+\frac{3}{2} }   } 
\right),
\end{equation}
where 
\[c=\frac{1}{\sqrt{2\pi}}
\int_{-\infty}^{\infty}|x|^3 e^{-\frac{1}{2}x^2} dx.\]
This together with (\ref{M26}) immediately imply
\begin{equation}\label{M30}
\mathbb{E} |Y_k|^3
\leq C\left(
\frac{1}{\lambda_k}\right)^{9/4},
\end{equation}
where $C$ is some positive constant. 

Recall that $\lambda_k=kF+\lambda$, where $\lambda \geq bN^2$ for some positive $b$ (Lemma \ref{L1}). Hence, bound (\ref{M30}) yields
\begin{equation}\label{M31}
\sum_{k=1}^N\mathbb{E} |Y_k|^3
\leq O(N^{-7/2}),
\end{equation}
which together with Corollary \ref{C1} implies the statement of the Proposition. \hfill$\Box$

To make use of the bound  (\ref{Pe}) we split the integral in  (\ref{M23}) into three parts:
\[
|f_{Z_N}(x)-\frac{1}{\sqrt{2\pi}}e^{-x^2/2}|\leq \frac{1}{\sqrt{2\pi} }
\int_{|t|\leq \frac{1}{4L_N}}|\phi_{Z_N}(t)  -e^{-t^2/2}|dt 
\]
\[
+\frac{1}{\sqrt{2\pi} } \int_{|t|>\frac{1}{4L_N}} |\phi_{Z_N}(t) | dt +
\frac{1}{\sqrt{2\pi} } \int_{|t|> \frac{1}{4L_N}} e^{-t^2/2} dt  .
\]
Bounds (\ref{Pe}) and (\ref{M32}) allow us to derive from here
\begin{equation}\label{M33}
|f_{Z_N}(x)-\frac{1}{\sqrt{2\pi}}e^{-x^2/2}|\leq  O(L_N) +\frac{1}{\sqrt{2\pi} } \int_{|t|>\frac{1}{4L_N}} |\phi_{Z_N}(t) | dt + O\left( e^{-1/L_N} \right)
\end{equation}
\[
=O(N^{-1/2}) +\frac{1}{\sqrt{2\pi} } \int_{|t|>\frac{1}{4L_N}} |\phi_{Z_N}(t) | dt.
\]

 Consider the remaining integral on the right in (\ref{M33}). Observe that by the definition (\ref{ZN})
\begin{equation}\label{M36}
|\phi_{Z_N}(t) | = \left|\prod _{k=1}^N\phi_{Y_k} \left( \frac{t}{\sqrt{\sigma^2_N}} \right) \right|= \prod _{k=1}^N\left|\phi_{X_{k, \lambda}}
\left( \frac{t}{\sqrt{\sigma^2_N}} \right) 
\right|.
\end{equation}
We shall derive now how fast  $|\phi_{X_{k, \lambda}}(t) |$ decays  in $|t|$.

By the Plancherel identity (see, e.g., \cite{F}) we have
\begin{equation}\label{M35}
\frac{1}{\sqrt{2\pi}} \int_{-\infty} ^{\infty}|\phi_{X_{k, \lambda}}(t) |^2 dt = 
\int_{-\infty} ^{\infty}f_{k, \lambda}(x)^2 dx,
\end{equation}
where by the definition (\ref{fXl})
\[
f_{k,\lambda}(u)=
\frac{e^{ 
-\frac{\beta }{u}  - (\lambda +F k) u 
}} 
{ \int_0^1  e^{-\lambda x -\frac{\beta}{x}} dx } 
, \ \ \ u  \in [0,1].
\]
Using again notation
\[\lambda_k=\lambda + kF,\]
we shall write here
\begin{equation}\label{fXk}
f_{k, \lambda}(u)=f_{\lambda_k}(u).
\end{equation}
Consider now for large $\lambda$
\[
 \int_{-\infty} ^{\infty}f_{\lambda}(x)^2 dx=\frac{\int_0^1 e^{-2\lambda x -
\frac{2\beta}{x}} dx}{
\left( \int_0^1  e^{-\lambda x -\frac{\beta}{x}} dx\right)^2
}=
\frac{\int_0^{\infty}  e^{-2\lambda x -
\frac{2\beta}{x}} dx +O(e^{-\lambda})}{
\left( \int_0^{\infty}  e^{-\lambda x -\frac{\beta}{x}} dx +O(e^{-\frac{1}{2}\lambda})\right)^2
}.
\]
Hence, using the integrals $I_{\alpha}(\lambda, \beta)$ defined in (\ref{Ia}), and then applying
the formula (\ref{M13})
 we get from here
\begin{equation}\label{M37}
 \int_{-\infty} ^{\infty}f_{\lambda}(x)^2 dx=
\frac{I_{1}(2\lambda, 2\beta)+O(e^{-\frac{1}{2}\lambda})}{(I_{1}(\lambda, \beta))^2
+O(e^{-\frac{1}{2}\lambda})}
\end{equation}

\[=
\frac{
\frac{2(2\beta)^{1/2} K_{1}(4\sqrt{\lambda\beta})}{(2\lambda)^{1/2}}
+O(e^{-\frac{1}{2}\lambda})}{\left(
\frac{2\beta^{1/2} K_{1}(2\sqrt{\lambda\beta})}{\lambda^{1/2}}
\right)^2
+O(e^{-\frac{1}{2}\lambda})}=\frac{\sqrt{\lambda}K_{1}(4\sqrt{\lambda\beta}) + O(e^{-\frac{1}{3}\lambda})
}{2\sqrt{\beta}\left( K_{1}(2\sqrt{\lambda\beta})
\right)^2
+ O(e^{-\frac{1}{3}\lambda})
}.
\]
Recall that 
the asymptotic  (\ref{M12}) yields
\[\frac{K_{1}(2z)}{(K_{1}(z))^2}=\sqrt{\frac{z}{\pi}} 
\left( 1+ O(z^{-1})\right), \]
which together with 
 (\ref{M37}) gives us
\begin{equation}\label{M38}
 \int_{-\infty} ^{\infty}f_{\lambda}(x)^2 dx=
\frac{\sqrt{\lambda}}{2\sqrt{\beta}} \ 
\sqrt{\frac{
2\sqrt{\lambda\beta}
}{\pi}} 
\left( 1+ O(\lambda^{-1/2})\right)=\frac{1}{\sqrt{2\pi} \ \beta^{1/4} } \ \lambda^{3/4}\left( 1+ O(\lambda^{-1/2})\right).
\end{equation}
Substituting this into 
 (\ref{M35}) we derive
\begin{equation}\label{M39}
\int_{-\infty} ^{\infty}|\phi_{X_{k, \lambda}}(t) |^2 dt = 
\beta^{-1/4} \ \lambda_k^{3/4}\left( 1+ O(\lambda_k^{-1/2})\right).
\end{equation}

Recall that here $\lambda_k=\lambda +kF$ where $F\leq 4\beta N$ and $\lambda$
is defined by Lemma \ref{L1}. Hence, there are positive constants $a<A$  such that 
\begin{equation}\label{M40}
 aN^2 < \lambda_k = \lambda +kF < A N^2.
\end{equation}
uniformly in $1\leq k \leq N$. Therefore we derive from (\ref{M39}) 
\begin{equation}\label{M41}
\int_{-\infty} ^{\infty}|\phi_{X_{k, \lambda}}(t) |^2 dt =  O(N^{3/2})
\end{equation}
uniformly in $1\leq k \leq N$, which yields 
\begin{equation}\label{M42}
|\phi_{X_{k, \lambda}}(t) |^2  =  o\left(\frac{N^{3/2}}{|t|} \right)
\end{equation}
 as $|t| \rightarrow \infty$, again uniformly in $k$. The last bound 
together with Proposition \ref{P2} (which tells us that $L_N =O(N^{-1/2})$)
implies 
\begin{equation}\label{M43}
\sup_{|t|>\frac{1}{4L_N}}\left|\phi_{X_{k, \lambda}}
\left( \frac{t}{\sqrt{\sigma^2_N}} \right) \right|^2 =
o\left( N^{3/2}
\sqrt{\sigma^2_N}N^{-1/2}
 \right)
\end{equation}
uniformly in $1\leq k \leq N$. By Corollary \ref{C1}  we have $\sigma^2_N=O(N^{-2})$. Therefore (\ref{M43}) yields
\begin{equation}\label{M44}
\sup_{|t|>\frac{1}{4L_N}}\left|\phi_{X_{k, \lambda}} 
\left( \frac{t}{\sqrt{\sigma^2_N}} \right) \right|^2 =
o\left( 1
 \right)
\end{equation}
uniformly in $1\leq k \leq N$. 

Finally, making use of bounds  (\ref{M41})  and (\ref{M43})
and taking into account (\ref{M36})
we can bound the second integral in   (\ref{M33}):
\begin{equation}\label{M45}
\frac{1}{\sqrt{2\pi} } \int_{|t|>\frac{1}{4 L_N}} |\phi_{Z_N}(t) | dt 
\leq \frac{1}{\sqrt{2\pi} } \int_{|t|>\frac{1}{4 L_N}}
\prod _{k=1}^N\left|\phi_{X_{k, \lambda}}
\left( \frac{t}{\sqrt{\sigma^2_N}} \right) 
\right|dt 
\end{equation}
\[\leq \frac{1}{\sqrt{2\pi} } \left( 
\max_{3\leq k \leq N}\sup_{|t|>\frac{1}{4L_N}}\left|\phi_{X_{k, \lambda}} 
\left( \frac{t}{\sqrt{\sigma^2_N}} \right) \right|
\right)^{N-2}
\int_{|t|>\frac{1}{4 L_N}}
\prod _{k=1}^2\left|\phi_{X_{k, \lambda}}
\left( \frac{t}{\sqrt{\sigma^2_N}} \right) 
\right|dt 
\]
\[\leq \gamma^{N}O(N^3),\]
for some positive $\gamma<1$.

Substituting  the last bound into  (\ref{M33}) 
we get the statement of Lemma \ref{L2}. \hfill$\Box$

\begin{cor}\label{C2}
Under assumptions of Lemma \ref{L2}  one has
\begin{equation}\label{M48}
f_{S_{N,\lambda}} (x)= \frac{1}{\sqrt{2\pi \sigma^2_N}} \ e^{-\frac{(x-1)^2}{2 \sigma^2_N}} + O(N^{-1/2})= \frac{1}{\sqrt{2\pi \sigma^2_N}} \ 
\left(e^{-\frac{(x-1)^2}{2 \sigma^2_N}} + O(N^{-3/2}) \right),
\end{equation}
as well as 
\begin{equation}\label{M49}
f_{S_{N,\hat{k}, \lambda}} (x)= \frac{1}{\sqrt{2\pi \sigma^2_N}} \left(
\ e^{-\frac{(x-1+\mathbb{E}X_{k,\lambda})^2}{2 \sigma^2_N}} + O(N^{-3/2})\right)
\end{equation}
uniformly in $1\leq k \leq N$ and $x\in \mathbb{R}.$
\end{cor}

\noindent 
{\bf Proof.} Formula (\ref{M48}) follows  immediately  by Lemma \ref{L2}
and formula (\ref{ZN}), while   (\ref{M49}) follows by the same argument and Corollary \ref{C1}.
\hfill$\Box$

\bigskip

Now we turn to the proof of the statements of Theorem \ref{T1}. 

\subsection*{Subcritical phase.}
Let us start with the statements  $(a)$ and $(b)$.
Applying  (\ref{cfX}) and using the result of the last Corollary we get 
\begin{equation}\label{cEX}
\mathbb{E} \{ X_k \mid \sum_{i=1}^{N} X_i =1 \} 
=  \int_{0}^{1} 
x  
\frac{
f_{k,\lambda}(x) f_{ S_{ N, \hat{k}, \lambda} }
 (1-x)
}{
f_{S_{N,\lambda}} (1)
}
dx
\end{equation}
\[
=  
\int_{0}^{1} 
x  f_{k,\lambda}(x) \left( e^{-\frac{(x-\mathbb{E}X_{k,\lambda})^2}{2 \sigma^2_N}} + O(N^{-3/2}) \right) dx\]
\[
= \int_{0}^{1} 
x  f_{k,\lambda}(x)  e^{-\frac{(x-\mathbb{E}X_{k,\lambda})^2}{2 \sigma^2_N}} dx + O(N^{-3/2}) \mathbb{E}X_{k,\lambda}.
\]
By Proposition \ref{P1} we have here 
\begin{equation}\label{JE}
\mathbb{E}X_{k,\lambda}  = \sqrt{\frac{\beta}{\lambda_k}}+O(1/\lambda_k),
\end{equation}
where by Lemma \ref{L1} 
\begin{equation}\label{M60}
\lambda_k=kF+\lambda
= \beta N^2  + \left( k-\frac{N}{2} \right) F+\beta \left( \frac{F}{4\beta }
\right)^2 +
O(N\log N),
\end{equation}
 and $\sigma^2_N=CN^{-2}(1+o(1))$ for some positive $C$ by Corollary \ref{C1}. 

Consider 
\begin{equation}\label{M50}
{\cal I}(\lambda_k) =  \int_{0}^{1} x
f_{k,\lambda}(x) e^{-\frac{(x-\mathbb{E}X_{k,\lambda})^2}{2 \sigma^2_N}} dx
=\frac{1}{I_1(\lambda_k,\beta)} \int_{0}^{1} x^{\alpha} e^{- \lambda_k x - \frac{\beta}{x}}
 e^{-\frac{\left(x-\sqrt{\frac{\beta}{\lambda_k}}+O(1/\lambda_k)\right)^2}{2 \sigma^2_N}} 
dx,
\end{equation}
where we used 
 notations (\ref{Ia}).
For all $\lambda >cN^2$ and any  positive constant $a$ we have here
\begin{equation}\label{M51}
{\cal I}(\lambda) = \frac{1}{I_1(\lambda,\beta)} \left(\int_{0}^{a/\sqrt{\lambda}}
 xe^{- \lambda x - \frac{\beta}{x}}
 e^{-\frac{\left(x-\sqrt{\frac{\beta}{\lambda}}+O(1/\lambda)\right)^2}{2 \sigma^2_N}} 
dx +O(e^{-a \sqrt{\lambda}/2})\right).
\end{equation}
Choosing now $a= 8\sqrt{\beta}$ we derive with a help of (\ref{M13}) and  (\ref{M12}) 
\begin{equation}\label{M52}
{\cal I}(\lambda) = \frac{1}{I_1(\lambda,\beta)} \int_{0}^{a/\sqrt{\lambda}}
 x e^{- \lambda x - \frac{\beta}{x} }
 e^{
-\frac{\left( x-\sqrt{ \frac{\beta}{\lambda}}\right)^2
}{2 \sigma^2_N}}
dx
\left( 1+O\left( \frac{1}{N} \right) 
\right)
+O(e^{-\sqrt{\beta\lambda}}).
\end{equation} 
Then applying nearly  same argument as in  (\ref{M6}) and (\ref{M7}), we get 
 from here 
 \begin{equation}\label{M53}
{\cal I}(\lambda) = \sqrt{ \frac{\beta}{\lambda} }\left(1+ O\left(\frac{\log N}{\sqrt{N}} \right)\right).
\end{equation}
This together with (\ref{M60}) 
yields 
\[
{\cal I}(\lambda_k) = \sqrt{ \frac{\beta}{\lambda_k} }\left(1+ O\left(\frac{\log N}{\sqrt{N}} \right)\right)
\]
\begin{equation}\label{M62}
=\sqrt{ \frac{1}{
 N^2 +  \left( k-\frac{N}{2}
\right)\frac{F}{\beta} + \frac{F^2}{16\beta^2 }
} }\left(1+ O\left(\frac{\log N}{\sqrt{N}} \right)\right)
\end{equation}
Hence, combining the last formula with  (\ref{cEX}) and (\ref{M50}) 
we get 
\begin{equation}\label{M61}
\mathbb{E}\{  X_k\mid \sum_{i=1}^{N} X_i =1 \}= \sqrt{ \frac{\beta}{\lambda_k} }\left(1+ O\left(\frac{\log N}{\sqrt{N}} \right)\right)= 
\mathbb{E}X_{k,\lambda}\left(1+ O\left(\frac{\log N}{\sqrt{N}} \right)\right)
\end{equation}
and then using the last formula in (\ref{M62})
it is straightforward
 to  derive the statements on the conditional expectations in $(a)$ and $(b)$ of Theorem \ref{T1}. 

In a similar to (\ref{cEX}) fashion consider 
\begin{equation}\label{cVX}
\mathbb{E} \left\{ \left(X_k - \sqrt{ \frac{\beta}{\lambda_k} }\right)^2 \mid \sum_{i=1}^{N} X_i =1 \right\} 
=  \int_{0}^{1} 
\left(x  - \sqrt{ \frac{\beta}{\lambda_k} }\right)^2 \ 
\frac{
f_{k,\lambda}(x) f_{ S_{ N, \hat{k}, \lambda} }
 (1-x)
}{
f_{S_{N,\lambda}} (1)
}
dx
\end{equation}
\[
=  
\int_{0}^{1} 
 \left(x  - \sqrt{ \frac{\beta}{\lambda_k} }\right)^2 \  f_{k,\lambda}(x)  e^{-\frac{(x-\mathbb{E}X_{k,\lambda})^2}{2 \sigma^2_N}} dx + O(N^{-3/2}) {Var}(X_{k,\lambda}),
\]
\[
 = \frac{1}{I_1(\lambda,\beta)} \int_{0}^{a/\sqrt{\lambda}}
 \left(x  - \sqrt{ \frac{\beta}{\lambda_k} }\right)^2  e^{- \lambda x - \frac{\beta}{x} }
 e^{
-\frac{\left( x-\sqrt{ \frac{\beta}{\lambda}}\right)^2
}{2 \sigma^2_N}}
dx
\left( 1+O\left( \frac{1}{N} \right) 
\right)
+O(N^{-3/2}) {Var}(X_{k,\lambda}), 
\]
with the same constant $a$ as in (\ref{M52}),
where  by (\ref{M16}) 
\begin{equation}\label{cV3}
Var (X_{k,\lambda})=\frac{\sqrt{\beta}}{2 \lambda_k^{3/2}}(1+o(1)).
\end{equation}
as $\lambda_k \rightarrow \infty$.
From (\ref{cVX}) applying again arguments as in  (\ref{M6}) and (\ref{M7}),
we derive
\begin{equation}\label{cVX2}
\mathbb{E} \left\{ \left(X_k - \sqrt{ \frac{\beta}{\lambda_k} }\right)^2 \mid \sum_{i=1}^{N} X_i =1 \right\} 
=  {Var}(X_{k,\lambda}) \left(1+ O\left(\frac{\log N}{\sqrt{N}} \right)\right).
\end{equation}
which together with  (\ref{cV3}) and  (\ref{M61}) 
(recall also Corollary \ref{C2}, part $(I)$ and (\ref{JE}))
yields both (\ref{TVa}) and  (\ref{TVb}).

\subsection*{Supercritical phase.}

In the case when $F/N=F_0\geq 4\beta $ we apply a different strategy (since the central limit theorem does not work here). 

Let $F/N=F_0$ be a constant such that $F_0>4\beta.$
Consider first  
\begin{equation}\label{cEX1}
\mathbb{E} \{ X_1 \mid \sum_{i=1}^{N} X_i =1 \} 
=  \int_{0}^{1} 
x  
\frac{
f_{1,\lambda}(x) f_{ (S_{ N}-X_{1, \lambda}), \lambda }
 (1-x)
}{
f_{S_{N,\lambda}} (1)
}
dx, 
\end{equation}
where by Lemma \ref{L1} $(d)$ 
\begin{equation}\label{cEX6}
f_{1,\lambda}(x) = \frac{e^{-\frac{\beta}{x} + \lambda_0 x}}{c}
\end{equation}
for some finite constant $\lambda_0>0$, and $c=c(\beta, \lambda_0)$.
Hence, here (see (\ref{M24}))
\begin{equation}\label{cEX2}
\mathbb{E} X_{1, \lambda} = 1-\sqrt{\frac{4 \beta}{F_0} }=: m_1 , \ \ Var(X_{1, \lambda})= : \sigma
\end{equation}
are some positive constants depending only on $\beta$ and $F_0$.
Let us also write here 
\[\Sigma_2= S_{ N}-X_{1, \lambda} = \sum_{k=2}^NX_{k, \lambda} .\]
By the Corollary \ref{C1}  $(III) $ the variance of $\Sigma_2$ decays as $N^{-3/2}$, therefore although we cannot apply central limit theorem some concentration results still hold.
Consider for $\alpha=0,1,2$ and any $\varepsilon<m_1 /2$
\begin{equation}\label{cEX3}
J(\alpha)
=  \int_{0}^{1} 
x^{\alpha}
f_{1,\lambda}(x) f_{ \Sigma_2}
 (1-x)
dx
\end{equation}
\[= \int_{|x-m_1|>\varepsilon}  O(1)  f_{ \Sigma_2}
 (1-x)
dx+
 (m_1+O(\varepsilon))^{\alpha}\int_{|x-m_1|\leq \varepsilon} f_{1,\lambda}(x) f_{ \Sigma_2}
 (1-x)dx.\]
Note that here by Corollary \ref{C1} (b) 
\begin{equation}\label{J8}
\int_{|x-m_1|>\varepsilon}   f_{ \Sigma_2}
 (1-x)
dx= \mathbb{P} \{ 
|\Sigma_2 -(1-m_1)|>\varepsilon
\} =O\left( \frac{1}{F^{3/2} \varepsilon^2}\right),
\end{equation}
while
\[\int_{|x-m_1|\leq \varepsilon} f_{1,\lambda}(x) f_{ \Sigma_2}
 (1-x)dx\geq \min_{|x-m_1|\leq \varepsilon} f_{1,\lambda}(x) \left(1-
O\left( \frac{1}{F^{3/2} \varepsilon^2}\right)\right)
.\]
Hence,
\begin{equation}\label{cEX5}
J(\alpha)
=  \int_{0}^{1} 
x^{\alpha}
f_{1,\lambda}(x) f_{ \Sigma_2}
 (1-x)
dx
\end{equation}
\[= \left( O\left( \frac{1}{F^{3/2} \varepsilon^2}\right)
 +
 (m_1+O(\varepsilon))^{\alpha}
\right)
\int_{0}^1 
f_{1,\lambda}(x) f_{ \Sigma_2}
 (1-x)dx \]
\[= \left( O\left( \frac{1}{F^{3/2} \varepsilon^2}\right)
 +
 (m_1+O(\varepsilon))^{\alpha}
\right) f_{S_{N,\lambda}} (1).\]
Making use of this formula in (\ref{cEX1}) we derive
\begin{equation}\label{cEX4}
\mathbb{E} \{ X_1 \mid \sum_{i=1}^{N} X_i =1 \} 
= m_1+O(\varepsilon) +  O\left( \frac{1}{F^{3/2} \varepsilon^2}\right).
\end{equation}
Choosing $\varepsilon=1 / \sqrt{N}$ we get from here 
\[\mathbb{E} \{ X_1 \mid \sum_{i=1}^{N} X_i =1 \} 
= m_1+O\left(\frac{1}{\sqrt{N}}\right),\]
which confirms (\ref{cEXsup}).

In a similar to (\ref{cEX4}) fashion we get as well for all positive $\varepsilon<m_1 /2$
\begin{equation}\label{cEX7}
\mathbb{E} \{ (X_1-m_1)^2 \mid \sum_{i=1}^{N} X_i =1 \} 
= O\left( \frac{1}{F^{3/2} \varepsilon^2}\right)
 + O(\varepsilon^{2}).
\end{equation}
Choosing this time $\varepsilon^2=1 / N^{3/4}$ we get 
\[\mathbb{E} \{ (X_1-m_1)^2\mid \sum_{i=1}^{N} X_i =1 \} 
= O\left(\frac{1}{N^{3/4}}\right),\]
which together with (\ref{cEXsup}) confirms (\ref{supV1}).

Next we consider 
$\mathbb{E} \{ X_k\mid \sum_{i=1}^{N} X_i =1 \}=
\{ X_{k,\lambda}\mid \sum_{i=1}^{N} X_{i,\lambda} =1 \} $ 
for $k>1$. Recall that here $\lambda_k=(k-1)F -  \lambda_0 $ 
where by the assumption $F=F_0 N $ and by the Lemma \ref{L1} $\lambda_0=\lambda_0(F_0)$ is also some constant.
The density of  $X_{k,\lambda}$ is
\[ f_{k,\lambda}(x) =\frac{e^{-\frac{\beta}{x} -\lambda_k x}}{
\int_0^1 e^ {-\frac{\beta}{x} -\lambda_k x } dx
}=:\frac{e^{-\frac{\beta}{x}  -(k-1)F  x+ \lambda_0 x } }{c_k}.\]
Let us denote here
\begin{equation}\label{lM}
 \ \mathbb{E} X_{k,\lambda}=m_k= 
\sqrt{
\frac{\beta}{(k-1)F+\lambda_0}
} +O\left(
\frac{1}{(k-1)F}
 \right),
\end{equation}
Let also 
\[\Sigma_{2, \hat{k}}= \Sigma_2 - X_{k,\lambda},\]
whose density we denote $f_{\Sigma_{2, \hat{k}}}$.
In these notations we have
\begin{equation}\label{cEX8}
\mathbb{E} \{ X_k\mid \sum_{i=1}^{N} X_i =1 \} 
= \frac{1}{
f_{S_{N,\lambda}} (1)
}
\int_{0}^{1} \left( \int_{0}^{z} 
x  
f_{k,\lambda}(x) f_{1, \lambda}(z-x)dx \right)\ f_{\Sigma_{2, \hat{k}}}(1-z)
dz.
\end{equation}
Consider now for $\alpha=0,1$
\[g_{\alpha}(z) = \int_{0}^{z} 
x^{\alpha}  
f_{k,\lambda}(x) f_{1, \lambda}(z-x)dx.\]
Using bound (\ref{M5}) we derive first for any $\varepsilon=o(1)$ and $z>m_k+\varepsilon/\sqrt{\lambda_k}$
\[
g_{\alpha}(z) = \int_{|x-m_k|>\varepsilon/\sqrt{\lambda_k}, \  0<x<z} 
x^{\alpha}
f_{k,\lambda}(x) f_{1, \lambda}(z-x)dx
+ \int_{|x-m_k|<\varepsilon/\sqrt{\lambda_k}} 
x^{\alpha}
f_{k,\lambda}(x) f_{1, \lambda}(z-x)dx
\]
\begin{equation}\label{JN}
= \left( 
O\left(e^{ -\frac{\sqrt{\lambda_k}}{4\sqrt{\beta}}\varepsilon^2 
}\right) + 
\left( m_k+O\left( \frac{\varepsilon}{\sqrt{\lambda_k}} \right) \right) ^{\alpha}\right) 
\int_{|x-m_k|<\varepsilon/\sqrt{\lambda_k}} 
f_{k,\lambda}(x) f_{1, \lambda}(z-x)dx,
\end{equation}
which yields both 
\begin{equation}\label{cEX9}
g_{\alpha}(z) =\left( 
O\left(e^{ -\frac{\sqrt{\lambda_k}}{4\sqrt{\beta}}\varepsilon^2 
}\right) + \left( 
m_k+O\left( \frac{\varepsilon}{\sqrt{\lambda_k}} \right) \right)^{\alpha}\right) 
g_{0}(z),
\end{equation}
and a bound (recall (\ref{cEX6}))
\begin{equation}\label{cEX11}
g_{\alpha}(z) \geq \left( 
O\left(e^{ -\frac{\sqrt{\lambda_k}}{4\sqrt{\beta}}\varepsilon^2 
}\right) + \left( 
m_k+O\left( \frac{\varepsilon}{\sqrt{\lambda_k}} \right)  \right)^{\alpha}\right)
\frac{
e^{-\frac{\beta}{z-(m_k+\varepsilon/\sqrt{\lambda_k})}}
}{c}
\end{equation}
\[\times \mathbb{P} \{ |X_{k,\lambda}-m_k|\leq 
\frac{\varepsilon}{\sqrt{\lambda_k}}
 \}. \]
Recall that here
\[Var(X_{k,\lambda})=\frac{\sqrt{\beta}}{2}
\frac{1}{( (k-1)F )^{3/2} }.\]
Choosing now 
\begin{equation}\label{vare}
\varepsilon= \frac{\log N }{(kF)^{1/4}}
\end{equation}
and  taking into account (\ref{lM}) 
we derive from (\ref{cEX11}) with a help of the Chebyshev's inequality
\begin{equation}\label{cEX12}
g_{\alpha}(z) \geq 
m_k^{\alpha}
\frac{
e^{-\frac{\beta}{z-(m_k+\varepsilon/\sqrt{\lambda_k})}}
}{2c} \left( 1-O\left( \frac{\lambda_k}{\varepsilon^2(kF)^{3/2}} \right) \right)
\geq 
m_k^{\alpha}
\frac{
e^{-\frac{\beta}{z-(m_k+\varepsilon/\sqrt{\lambda_k})}}
}{4c} .
\end{equation}
On the other hand, for all $z$ it holds that
\[g_{1} (z) = O(z^2).\]
With a help of  the last bound and (\ref{cEX9}) consider the integral in  (\ref{cEX8})
with the same choice of $\varepsilon$ as in (\ref{vare})
\begin{equation}\label{cEX10}
\int_{0}^{1} g_1(z)\ f_{\Sigma_{2, \hat{k}}}(1-z)
dz= 
\int_{0}^{m_k+\varepsilon/\sqrt{\lambda_k}} 
O(z^2)
\ f_{\Sigma_{2, \hat{k}}}(1-z)
dz
\end{equation}
\[+\int_{m_k+\varepsilon/\sqrt{\lambda_k}}^1
\left( 
O\left(e^{ -\frac{\sqrt{\lambda_k}}{4\sqrt{\beta}}\varepsilon^2 
}\right) + 
m_k+O\left( \frac{\varepsilon}{\sqrt{\lambda_k}} \right) \right) 
g_{0}(z)
\ f_{\Sigma_{2, \hat{k}}}(1-z)
dz 
\]
\[= O(m_k^2) + \left(m_k+O\left( \frac{\varepsilon}{\sqrt{\lambda_k}} \right) \right) 
\int_{m_k+\varepsilon/\sqrt{\lambda_k}}^1g_{0}(z)
\ f_{\Sigma_{2, \hat{k}}}(1-z)
dz,
\]
where by   (\ref{cEX12})  we have
\[\int_{m_k+\varepsilon/\sqrt{\lambda_k} }^1 
g_{0}(z)
\ f_{\Sigma_{2, \hat{k} } }(1-z)
dz\]
\[
\geq 
\int_{m_k+m_1 /2 }^1 
g_{0}(z)
\ f_{\Sigma_{2, \hat{k} } }(1-z)
dz
\geq a
\mathbb{P} 
\{ 
| \Sigma_{2, \hat{k} }-(1-m_k)|\geq m_1 /2
 \} \geq b \]
for some positive constants $a$ and $b$ uniformly in $N$. 
The last bound yields in turn together with (\ref{cEX10})
\begin{equation}\label{cEX14}
\int_{0}^{1} g_1(z)\ f_{\Sigma_{2, \hat{k}}}(1-z)
dz= \left(
m_k+O
\left( \frac{\log N}{ \lambda_k^{3/4} } \right) 
\right) 
\int_{0}^1g_{0}(z)
\ f_{\Sigma_{2, \hat{k}}}(1-z)
dz 
\end{equation}
\[
= \left(
m_k+O
\left( \frac{\log N}{ \lambda_k^{3/4} } \right) 
\right) f_{S_{N,\lambda}}(1).
\]
Substituting the last formula into (\ref{cEX8}) we get 
\begin{equation}\label{cEX13}
\mathbb{E} \{ X_k\mid \sum_{i=1}^{N} X_i =1 \} 
= 
m_k+O
\left( \frac{\log N}{ \lambda_k^{3/4} } \right).
\end{equation}
This together with (\ref{lM}) proves (\ref{supE}).

Consider now
\begin{equation}\label{cEX15}
\mathbb{E} \{ (X_k-m_k)^2\mid \sum_{i=1}^{N} X_i =1 \} 
\end{equation}
\[= \frac{1}{
f_{S_{N,\lambda}} (1)
}
\int_{0}^{1} \left( \int_{0}^{z} 
(x-  m_k)^2
f_{k,\lambda}(x) f_{1, \lambda}(z-x)dx \right)\ f_{\Sigma_{2, \hat{k}}}(1-z)
dz.
\]
Using again the same argument as above first we derive for all $z>2 m_k$,  $\varepsilon=o(1)$ but such that $\sqrt{\lambda_k}\varepsilon^2 \gg 1$, and any $\alpha =0,1$
\[
 G_{\alpha }(z):=\int_{0}^{z} 
(x-  m_k)^{\alpha }
f_{k,\lambda}(x) f_{1, \lambda}(z-x)dx = 
\]
\[
= \left( 
O\left(e^{ -\frac{\sqrt{\lambda_k}}{4\sqrt{\beta}}\varepsilon^2 
}\right) +1\right) 
\int_{|x-m_k|<\varepsilon/\sqrt{\lambda_k}} (x-  m_k)^{\alpha }
f_{k,\lambda}(x) f_{1, \lambda}(z-x)dx
\]
Then similar to the derivation of (\ref{M9}) we get from here
\begin{equation}\label{cEX16}
 G_{\alpha }(z)
= \left( 
O\left(e^{ -\frac{\sqrt{\lambda_k}}{4\sqrt{\beta}}\varepsilon^2 
}\right) +1\right) 
O(Var (X_{k,\lambda}))^{\alpha/2 }
 f_{1, \lambda}(z-m_k)(1  +O(\varepsilon/\sqrt{\lambda_k})),
\end{equation}
where we used the fact that  in the density $f_{1, \lambda}$ (see (\ref{cEX6})) parameter $\lambda_0$
is $O(1)$.
The last formula allows us to derive from (\ref{cEX15})
\begin{equation}\label{J11}
\mathbb{E} \{ (X_k-m_k)^2\mid \sum_{i=1}^{N} X_i =1 \} \leq 
\frac{(1+o(1))
\int_{2m_k}^{1} G_{2 }(z)\ f_{\Sigma_{2, \hat{k}}}(1-z)
dz}{\int_{0}^{1} G_{0 }(z)\ f_{\Sigma_{2, \hat{k}}}(1-z)dz}
\end{equation}
\[=O(Var (X_{k,\lambda}))
=O\left( \frac{1}{(kF)^{3/2}} \right),
\]
which together with (\ref{supE}) yields (\ref{supV}).

Let $F/N=4\beta.$ Then similar to the previous case we consider 
\begin{equation}\label{J6}
\mathbb{E} \{ X_k \mid \sum_{i=1}^{N} X_k =1 \} 
=  \int_{0}^{1} 
x  \ 
\frac{
f_{k,\lambda}(x) f_{  S_{ N, \hat{k}, \lambda }}
 (1-x)
}{
f_{S_{N,\lambda}} (1)
}
dx, 
\end{equation}
where $\lambda$ is chosen according to Lemma \ref {L1}, so that
\begin{equation}\label{J4}
f_{k,\lambda}(x) = \frac{e^{-\frac{\beta}{x} - \lambda_k x}}{c_k}
\end{equation}
with $\lambda_k = kF+\Theta(N)=4\beta Nk+\Theta(N),$ positive. Thus by Lemma \ref{L1} and Proposition \ref{P1} we have here
\begin{equation}\label{J13}
m_k= EX_{k,\lambda}= \sqrt{\frac{\beta}{\lambda_k }}(1+o(1))= \sqrt{\frac{1}{4 Nk+\Theta(N) }}(1+o(1)).
\end{equation}

Similar to (\ref{cEX3}) consider first for any $\varepsilon>0$
\begin{equation}\label{J9}
J
=  \int_{0}^{1} 
x
f_{k,\lambda}(x) f_{ S_{ N, \hat{k}, \lambda }}
 (1-x)
dx
\end{equation}
\[=  \int_{|x-m_k|>\varepsilon}x  f_{k,\lambda}(x)  f_{S_{ N, \hat{k}, \lambda }  }
 (1-x)
dx+
\int_{|x-m_k|\leq \varepsilon} x f_{k,\lambda}(x) 
f_{ S_{ N, \hat{k}, \lambda }}
 (1-x)dx.\]
Then with a help of Corollary \ref{C1} (part (II)) we derive as in (\ref{J8}) (and using also the Chebyshev's inequality)
\[
\int_{|x-m_k|>\varepsilon} x f_{k,\lambda}(x)  f_{S_{ N, \hat{k}, \lambda }  }
 (1-x)
dx \leq 
\max_{|x-m_k|\geq \varepsilon} (x f_{k,\lambda}(x) )
\mathbb{P} \{ 
| S_{ N, \hat{k}, \lambda } -(1-m_k)|>\varepsilon
\} 
\]
\begin{equation}\label{J10}
= \max_{|x-m_k|\geq \varepsilon} (x f_{k,\lambda}(x)) \ O\left( \frac{1}{F^{3/2} \varepsilon^2}\right),
\end{equation}
while
\[\int_{|x-m_1|\leq \varepsilon} x f_{1,\lambda}(x) f_{ S_{ N, \hat{k}, \lambda }}
 (1-x)dx \geq \min_{|x-m_k|\leq \varepsilon} (x f_{k,\lambda}(x)) \left(1-
O\left( \frac{1}{F^{3/2} \varepsilon^2}\right)\right)
.\]
Using the exact form (\ref{J4}), one can see that for any $\varepsilon \geq 1/\lambda_k$ we have here
\[ \min_{|x-m_k|\leq \varepsilon}(x f_{k,\lambda}(x) )=\max_{|x-m_k|\geq \varepsilon}(x f_{k,\lambda}(x)). \]
This allows us to derive from  (\ref{J9}) similar to (\ref{cEX5})
\begin{equation}\label{J12}
J
= \left(  O\left( \frac{1}{F^{3/2} \varepsilon^2}\right)
 +1\right) \int_{|x-m_1|\leq \varepsilon} x f_{1,\lambda}(x) f_{ S_{ N, \hat{k}, \lambda }}
 (1-x)dx 
\end{equation}
\[= \left(  O\left( \frac{1}{F^{3/2} \varepsilon^2}\right)
 +1\right)
 (m_k+O(\varepsilon)) f_{ S_{ N,  \lambda }}
 (1).\]
Setting here  $\varepsilon = F^{-2/3} \geq 1/\lambda_k$ for all $k$ and large $N$,
we get 
\[J
= (m_k+O(F^{-2/3})) f_{ S_{ N,  \lambda }}
 (1).\]
Substituting this result into
 (\ref{J6}) and taking into account (\ref{J13}) we get (\ref{crEk}).

Next under the same condition $F/N=4\beta$ consider 
\begin{equation}\label{J14}
\mathbb{E} \{ (X_k-m_k)^2 \mid \sum_{i=1}^{N} X_k =1 \} 
=  \int_{0}^{1} 
(x-m_k) ^2 \ 
\frac{
f_{k,\lambda}(x) f_{  S_{ N, \hat{k}, \lambda }}
 (1-x)
}{
f_{S_{N,\lambda}} (1)
}
dx.
\end{equation}
With exactly same argument as we derived (\ref{J12})  we get here
\begin{equation}\label{J15}
 \int_{0}^{1} 
(x-m_k) ^2 \ 
f_{k,\lambda}(x) f_{  S_{ N, \hat{k}, \lambda }}
 (1-x)
dx= \left(  O\left( \frac{1}{N^{3/2} \varepsilon^2}\right)
 +1\right)
 O(\varepsilon^2) f_{ S_{ N,  \lambda }}
 (1).
\end{equation}
Setting  $\varepsilon = N^{-3/4}$ (which is also greater than  $\lambda_k$ for all $k$ and large $N$) and using the  result in
 (\ref{J14}) we obtain 
\[\mathbb{E} \{ (X_k-m_k)^2 \mid \sum_{i=1}^{N} X_k =1 \} 
= O(N^{-3/2}).\]
This together with (\ref{crEk}) yields (\ref{TVc}).

Finally, assume that  $F\gg N$. Similar to the case when $F=4\beta N$ we consider first 
\begin{equation}\label{J16}
\mathbb{E} \{ X_1 \mid \sum_{i=1}^{N} X_k =1 \} 
=  \int_{0}^{1} 
x \ 
\frac{
f_{1,\lambda}(x) f_{ \Sigma_2 }
 (1-x)
}{
f_{S_{N,\lambda}} (1)
}
dx,
\end{equation}
where $\Sigma_2:= S_{ N, \lambda }-X_{1, \lambda}$. In (\ref{J16})  we have 
\[f_{1,\lambda}(x) = \frac{e^{\lambda_1 x - \frac{\beta}{x}}}{c} \]
where (see Remark \ref{R1}) $\lambda_1 = \Theta(\sqrt{F/N})$, and $c$ is a normalizing constant. Recall also that by Lemma \ref{L1} $(d)$
\begin{equation}\label{J17}
\mathbb{E}X_{1, \lambda} = 1-\sqrt{\frac{4 \beta N}{F} },
\end{equation}
and by Corollary \ref{C2} part (III) 
\[Var (\Sigma_2  ) =O\left( \frac{1}{ F^{3/2}}\right).\]
Therefore
the same method as we used in the previous case works here exactly same way. 
Theorem is proved. \hfill$\Box$

\bigskip

{\bf Acknowledgments.} The author thanks V. Malyshev for the 
introducing this problem, and Y. Ameur for the helpful and inspiring
discussions on the general framework of this model.
The author also  gratefully acknowledges  discussions with D. Mason and  
B. S{\"o}derberg.

\end{document}